\documentclass[12 pt]{report}
\usepackage{amssymb,latexsym,epsfig,mathrsfs,graphpap,graphics,amsmath,amssymb}
\usepackage{amsmath}

\begin{document}

\vspace{.2in}\parindent=0mm

\begin{flushleft}
{\bf\Large { Fractional Biorthogonal wavelets in $L^2(\mathbb R)$ }}

  \parindent=0mm \vspace{.3in}
{\bf{  Owais Ahmad$^{*}$, N. A. Sheikh$^{1}$ and Firdous A. Shah$^{2}$}}
\end{flushleft}

\parindent=0mm \vspace{.1in}
{{\it $^{*}$Department of  Mathematics,  National Institute of Technology, Hazratbal, Srinagar -190 006, Jammu and Kashmir, India. E-mail: $\text{siawoahmad@gmail.com}$}}

\parindent=0mm \vspace{.1in}
{{\it $^{1}$Department of  Mathematics,  National Institute of Technology, Hazratbal, Srinagar -190 006, Jammu and Kashmir, India.E-mail: $\text{neyaznit@yahoo.co.in}$}}

\parindent=0mm \vspace{.1in}
{{\it\small$^{2}${Department of  Mathematics,  University of Kashmir, South Campus, Anantnag-192 101, Jammu and Kashmir, India. E-mail: $\text{fashah79@gmail.com}$}}

\parindent=0mm \vspace{.2in}
{\bf{Abstract.}}The fractional Fourier transform (FrFT), which is a generalization of the Fourier transform, has become the focus of many research papers in recent years because of its applications in electrical engineering and optics. In this paper, we introduce the notion of  fractional biorthogonal wavelets on $\mathbb{R}$ and obtain the necessary and sufficient conditions for the translates of a single function to form the fractional  Riesz bases for their closed linear span. We also provide a complete characterization for the  fractional biorthogonality of the translates of fractional scaling functions of two fractional MRA’s and the associated fractional biorthogonal wavelet families. Moreover, under mild assumptions on the fractional scaling functions and the corresponding fractional wavelets, we show that the fractional wavelets can generate Reisz bases for $L^2(\mathbb R).$

\parindent=0mm \vspace{.2in}
{\bf{Keywords:}}}  Frame; Fractional Biorthogonal  wavelets; Fractional MRA; Fractional Fourier transform.

\parindent=0mm \vspace{.1in}


\parindent=0mm \vspace{.1in}
{\bf{1. Introduction}}

\parindent=0mm \vspace{.1in}
The Fourier transform has been used for more than a century in a wide range of applications. However, more recently, it was shown that the Fourier transform is inadequate for describing some physical applications or dealing with their underlying mathematical problems. As a result, some off-shoots of the Fourier transform, such as the windowed Fourier transform, the wavelet transform, and the fractional Fourier transform (FrFT) have been introduced to address the shortcoming of the Fourier transform. The FrFT, which is a generalization of the Fourier transform,  has gained considerable attention in the last 20 years or so because of its important applications in signal analysis, optics, and signal recovery and also because of its ability to treat some mathematical problems that could not otherwise be handled by the standard Fourier transform \cite{7}. The FrFT appeared implicitly in the work of N. Wiener in 1929 \cite{wn} as a way to solve certain types of ordinary and partial differential equations arising in quantum mechanics. Unaware of Wiener’s work, V. Namias in 1980 \cite{10} introduced the transform, which he called the FrFT, also to solve ordinary and partial differential equations arising in quantum mechanics from classical quadratic Hamiltonians. His work was later refined by McBride and Kerr \cite{mk}.Besides lot of advantages, the FrFT has one major drawback due to using global kernel i.e., the fractional Fourier representation only provides such  FrFT spectral content with no indication about the time localization of the FrFT spectral components. On the other hand, the  short-time FrFT has rectified almost all the limitations of FrFT, still in some cases short-time FrFT is also not applicable as in the case of real signals having  high spectral components for short durations and low spectral components for long durations. Therefore,  in order to obtain joint signal representations in both time and FrFT domains, Mendlovic et al. \cite{ men} first introduced the fractional wavelet transform (FrWT) in the context of time-frequency analysis. The  FrWT inherits the excellent mathematical properties of wavelet transform and FrFT along with some fascinating properties of its own. The idea behind this transform is deriving the fractional spectrum of the signal by using the FrFT and performing the wavelet transform of the fractional spectrum. Besides being a generalization of the wavelet transform, the FrWT can be interpreted as a rotation of the time–frequency plane and has been proved to relate to other time-varying signal analysis tools, which make it as a unified time–frequency transform. In recent years, this transform has been paid a considerable amount of attention, resulting in many applications in the areas of optics, quantum mechanics, pattern recognition and signal processing. For more about fractional wavelet transforms and their applications to signal and image processing, we refer to  \cite{{ofwpf},{dzw},{men},{ap}}.

\parindent=8mm \vspace{.1in}
 Along with the study of wavelet transforms, there had been a continuing research effort in the study of biorthogonal wavelets and their promising features in applications have attracted a great deal of interest  in recent years to extensively study them. During the late 1990’s, biorthogonal wavelets brought a major breakthrough into image compression, thanks to their natural feature of concentrating energy in a few transform coefficients. In traditional wavelet theory, biorthogonal wavelets have many advantages over orthogonal wavelets, by relaxing orthonormal to biorthogonal, additional degrees of freedom are added to design problems. Biorthogonal wavelets in $L^2(\mathbb R)$ were investigated by Bownik and Garrigos \cite{1}, Cohen et al. \cite{2}, Chui and Wang \cite{3} and many others. 
 
 \parindent=8mm\vspace{.1in}
Although there are many results for biorthogonal wavelets on the real-line $\mathbb R$, the counterparts on the fractional case are not reported yet in the literature. So this paper is concerned with the construction of fractional biorthogonal wavelets on $\mathbb R$. We introduce the notion of  fractional biorthogonal wavelets on $\mathbb{R}$ and obtain the necessary and sufficient conditions for the translates of a single function to form the fractional  Riesz bases for their closed linear span. We also provide a complete characterization for the  fractional biorthogonality of the translates of fractional scaling functions of two fractional MRA’s and the associated fractional biorthogonal wavelet families. Moreover, under mild assumptions on the fractional scaling functions and the corresponding fractional wavelets, we show that the fractional wavelets can generate Reisz bases for $L^2(\mathbb R).$

 \parindent=8mm\vspace{.1in}
The article is structured in the following manner. In Section 2, we recall the basic definitions of fractional Fourier transform and fractional wavelet transform. In Section 3, we establish necessary and sufficient conditions for the translates of a function to form a fractional Riesz basis for its closed linear span.  In section 4, we give the definition of a fractional  MRA. We also define the projection operators associated with the fractional MRAs and show that they are uniformly bounded on $L^2(\mathbb R)$. In the concluding Section, we show that the fractional wavelets associated with fractional dual MRA’s are biorthogonal and generate Riesz bases for $L^2(\mathbb R).$

\newpage
\parindent=0mm \vspace{.2in}
{\bf{2. Fractional Fourier and wavelet transforms}}

\parindent=0mm \vspace{.1in}
This section gives the basic background to the theory of fractional Fourier and wavelet transforms which is as follows.

\parindent=8mm \vspace{.1in}

In 1980, Victor Namias \cite{10} introduced the concept of fractional Fourier transform (FrFT) as a generalization of the conventional Fourier transform to  solve certain problems arising in quantum mechanics.  It is also referred as {\it rotational Fourier transform} or {\it angular Fourier transform} since it depends on a parameter $\alpha$ which is interpreted as a rotation by an angle $\alpha$ in the time-frequency plane. Like the ordinary Fourier transform  corresponds to a rotation in the time frequency plane over an angle $\alpha = 1\times {\pi}/{2}$, the FrFT corresponds to a rotation over an arbitrary angle $\alpha= \rho\times {\pi}/{2}$ with $\rho\in\mathbb R$.

\parindent=0mm\vspace{.2in}
The fractional Fourier transform with parameter $\alpha$  of function $f(t)$ is defined by

$$\mathcal{F}_\alpha\big\{f(t)\big\}(\xi)=\hat f^\alpha(\xi)=\int_{-\infty}^{\infty} {\mathcal K}_\alpha (t,\xi)f(t)\,dt,\eqno(2.1)$$

\parindent=0mm \vspace{.1in}
where $  {\cal K}_\alpha (t,\xi)$ is the so-called kernel of the FrFT given by
$${\cal K}_\alpha (t,\xi)=\left\{\begin{array}{ll}
  C_\alpha \exp\Big\{i(t^2+\xi^2)\dfrac{\cot\alpha}{2}-it\xi \csc\,\alpha\Big\},&\alpha\neq n\pi, \\
\delta(t-\xi), &\alpha=2n\pi, \\
\delta(t+\xi), &\alpha=(2n\pm 1)\pi, \\
\end{array}{}\right.\eqno(2.2)$$

\parindent=0mm \vspace{.1in}
$\alpha={\rho\pi}/{2}$ denotes the rotation angle of the transformed signal for FrFT,   the FrFT operator is designated by $\mathcal{F}_\alpha$ and
\begin{align*}
C_\alpha=\left(2\pi i \sin\alpha\right)^{-1/2}e^{i\alpha/2}=\sqrt{\dfrac{1-i\cot\alpha}{2\pi}}.\tag{2.3}
\end{align*}

\parindent=0mm\vspace{.0in}
The corresponding inversion formula is given by
$$f(t)=\int_{-\infty}^{\infty}\overline{{\cal K}_{\alpha}(t,\xi)}\,\hat{f}^\alpha(\xi)\,d\xi,\eqno(2.4)$$
where
\begin{align*}
  \nonumber {\cal K}_{\alpha}(t,\xi) &= \frac{(2\pi i\sin\alpha)^{1/2}\,e^{-i\alpha/2}}{\sin\alpha}\cdot \exp\left\{\frac{-i(t^2+\xi^2)\cot\alpha}{2}+it\xi \csc\,\alpha\right\}\\
  \nonumber &=\overline{C_\alpha}\exp\left\{\frac{-i(t^2+\xi^2)\cot\alpha}{2}+i t\xi \csc\,\alpha\right\}\\
   &={\cal K}_{-\alpha}(t,\xi)\tag{2.5}
\end{align*}
and
$$C_{\alpha}=\frac{(2\pi i\sin\alpha)^{1/2}e^{-i\alpha/2}}{2\pi \sin\alpha}=\sqrt{\dfrac{1+i\cot\alpha}{2\pi}}=C_{-\alpha}.\eqno(2.6)$$

\parindent=0mm \vspace{.0in}
{\bf {Definition 2.1.}} A  fractional wavelet is a function $\psi\in L^2(\mathbb{R})$ which satisfies the following condition:
$$C_{\psi}^{\alpha}=\int_{\mathbb{R}}\frac{\left|\mathcal{F}_\alpha\left\{e^{{-i(t-\xi)^2}/{2}\,\cot\alpha}\psi\right\}(\xi) \right|^2}{|\xi|}\,d\xi<\infty,\eqno(2.7)$$

\parindent=0mm\vspace{.1in}
where $\mathcal{F}_\alpha$ denotes the FrFT operator.

\parindent=8mm\vspace{.1in}
Analogous to the classical wavelets, the fractional wavelets can be obtained from a fractional mother wavelet $\psi\in L^2(\mathbb{R})$ by the combined action of translation and dilations as
$$\psi_{{\alpha},a,b}(t)=\frac{1}{\sqrt{a}}\,\psi\left(\frac{t-b}{a}\right)\exp\left\{\frac{-i(t^2-b^2)\cot\alpha}{2}\right\}\eqno(2.8)$$

\parindent=0mm\vspace{.1in}
where $a\in\mathbb R^+$ and $b\in\mathbb R$ are scaling and translation parameters, respectively. 
If $\alpha={\pi}/{2}$, then $\psi_{\alpha, a,b}$ reduces to the conventional wavelet basis.

\parindent=5mm\vspace{.1in}
Note that if $\psi(t)\in L^2(\mathbb{R})$, then $\psi_{\alpha, a,b}(t)\in L^2(\mathbb{R})$,
$$\left\|\psi_{\alpha, a,b}\right\|^2_{2}=|a|^{-1}\int_{-\infty}^{\infty}\left|\psi\left( \frac{t-b}{a}\right)\right|^2dt = \int_{-\infty}^{\infty}\big|\psi(y)\big|^2dy =\big\|\psi\big\|_{2}^{2}.$$

\parindent=0mm\vspace{.1in}
Moreover, the fractional Fourier transform of $\psi_{\alpha, a,b}(t)$ is given by

\parindent=0mm\vspace{.1in}
$\mathcal{F}_\alpha\big\{\psi_{\alpha, a,b}(t)\big\}$
$$=\sqrt{a}\,\exp\left\{\frac{i(b^2+\xi^2)\cot\alpha}{2}-ib\,\xi \csc\,\alpha-\frac{ia^2\xi^2\cot\alpha}{2}\right\}\mathcal{F}_\alpha\Big\{ e^{{-i(\cdot)^2\cot\alpha}/{2}}\psi\Big\}(a\xi)\eqno(2.9)$$

\parindent=5mm\vspace{.1in}
The continuous fractional wavelet transform (FrWT) of  function $f\in L^2(\mathbb R)$ with respect to an analyzing wavelet $\psi\in L^2(\mathbb R)$ is defined as
$${\mathscr W}^{\alpha}_{\psi}f(a,b)=\big\langle f,\psi_{\alpha, a,b}\big\rangle= \frac{1}{\sqrt{a}}\int_{-\infty}^{\infty}f(t)\,\overline{\psi\left( \frac{t-b}{a}\right)} \exp\left\{ \frac{i(t^2-b^2)\cot\alpha}{2}\right\} dt\eqno(2.10)$$

\parindent=0mm\vspace{.1in}
where $\psi_{\alpha,a,b}(t)\in L^2(\mathbb R)$ is given by (2.8).

\parindent=5mm\vspace{.2in}
The FrWT (2.10) deals generally with continuous functions, i.e. functions which are defined at all values of the time $t$. However, in many applications, especially in signal processing, data are represented by a finite number of values, so it is important and often useful to consider the discrete version of the continuous FrWT  (2.10). From a mathematical point of view, the continuous parameters $a$ and $b$ in (2.8) can be converted into a discrete one by assuming that $a$ and $b$ take only integral values. For a good discritization of the  wavelets, we choose $a=a_{0}^{-j}$ and $b=kb_{0}a_{0}^{-j}$, where $a_{0}$ and $b_{0}$  are fixed positive constants. Hence, the discritized wavelet family is defined as

$$\psi_{\alpha,j,k}(t)=a_{0}^{j/2}\,\psi\left(a_{0}^{j}t-kb_{0}\right)\exp\left\{-i~\frac{t^2-\big(kb_{0}a_{0}^{-j}\big)^2}{2}\cot\alpha\right\} \eqno(2.11)$$

\parindent=0mm\vspace{.1in}
where the integers $j$ and $k$ are the controlling factors for the dilation and translation, respectively and are contained in a set of integers. For computational efficiency, the discrete wavelet parameters $a_{0}=2$ and $b_{0}=1$ are commonly used so that equation (2.11) becomes
$${\mathscr F}_{\psi}(j,k):=\left\{\psi_{\alpha, j,k}(t)=2^{j/2}\,\psi\left(2^{j}t-k\right)e^{-i~\frac{t^2-(k2^{-j})^2}{2}\cot\alpha}, j,k\in\mathbb Z\right\} . \eqno(2.12)$$

\parindent=0mm\vspace{.1in}
The fractional wavelet system ${\mathscr F}_{\psi}(j,k)$ is called a {\it fractional wavelet frame}, if there exist positive constants $A$ and
$B$ such that

$$A\big\|f \big\|^2_{2} \le \sum_{j\in\mathbb Z}\sum_{k\in \mathbb Z} \left|\big\langle f, \psi_{\alpha, j, k}\big\rangle\right|^2 \le B \big\|f\big\|^2_{2},\eqno(2.13)$$

\parindent=0mm \vspace{.1in}
holds for every $f\in  L^2(\mathbb R)$, and we call the optimal constants $A$ and $B$ the lower frame bound and the upper frame bound, respectively. A {\it tight fractional wavelet frame} refers to the case when $A = B$, and a Parseval  frame refers to the case when $A = B = 1$. On the other hand if only the right hand side of the above double inequality holds, then we say ${\mathscr F}_{\psi}(j,k)$ a {\it Bessel system}.

\parindent=0mm \vspace{.1in}
{\bf{3. Reisz Bases of Translates}} 

\parindent=0mm \vspace{.1in}
\textbf{Definition 3.1.}  Let $\{\psi_m :m \in \mathbb Z\}$ and $\{\widetilde{\psi}_n : n\in \mathbb Z\}$ be two collections of functions in $L^2(\mathbb R)$. We say that they are orthogonal if

$$\left\langle \psi_n, \widetilde{\psi}_m \right\rangle = \delta_{n,m}~~\forall~ m,n \in \mathbb Z.$$

\parindent=0mm \vspace{.0in}
\textbf{Definition 3.2.} A collection of functions $\{\psi_n :n \in \mathbb Z\}$ in $L^2(\mathbb R)$ is said to be linearly independent if there exists a coefficient sequence $a[n] \in \ell^2(\mathbb Z)$ such that
$$\sum_{n=1}^{\infty} a[n]\psi_n = 0~\textit{in}~L^2(\mathbb R),$$
then $a[n] =0~~\forall n \in \mathbb N.$

\parindent=0mm \vspace{.1in}
\textbf{Lemma 3.3.} Let $\{\psi_n :n \in \mathbb Z\}$ be a collection of functions in $L^2(\mathbb R)$. Suppose that there is a collection $\{\widetilde{\psi}_n : n \in \mathbb Z\}$ in $L^2(\mathbb R)$ which is orthogonal to $\{\psi_n : n\in \mathbb Z\}$.Then $\{\psi_n : n\in \mathbb Z\}$ is linearly independent.

\parindent=0mm \vspace{.1in}
\textbf{Proof.} Let $a[n] \in \ell^2(\mathbb Z)$ be a coefficient sequence satisfying
$$\sum_{n=1}^{\infty} a[n] \psi_n  = 0 ~\textit{in}~L^2(\mathbb R).$$
Then for each $m \in \mathbb N,$ we have
\begin{align*}
0 &= \langle 0 , \widetilde{\psi}_m\rangle\\\
& = \big\langle \sum_{n=1}^{\infty} a[n] \psi_n, \widetilde{\psi}_n\big\rangle \\\
& =\sum_{n=1}^{\infty} a[n]\langle \psi_n, \widetilde{\psi}_n\rangle \\
& =a[m].
\end{align*}
Hence $\{\psi_n : n\in \mathbb Z\}$ is linearly independent.

\parindent=0mm \vspace{.1in}
\textbf{Definition 3.4.} A collection of functions $\{g_n(x)\}$ in $L^2(\mathbb R)$ is said to form a Reisz basis for a Hilbert space $\cal{H}$ if

\parindent=0mm \vspace{.1in}
(a) $\{g_n(x)\}$ is linearly independent, and

\parindent=0mm \vspace{.1in}
(b) there exists positive constants $A, B $ such that

$$A\|f\|_2^2 \le \sum_{n=1}^{\infty}\left|\langle f, g_n\rangle\right|^2 \le B\|f\|_2^2~~\forall~~ f \in \cal{H}.$$

\parindent=0mm \vspace{.1in}
In the following lemma, we establish a necessary and sufficient condition for the trans-
lates of two functions to be biorthogonal in fractional sense.

\parindent=0mm \vspace{.1in}
\textbf{Lemma 3.5.} Let the functions $\phi_\alpha(t-n)$ and $\widetilde{\phi}_\alpha (t-n)$ in $L^2(\mathbb R)$ are given. Then $\{\phi_\alpha (t-n) : n \in \mathbb Z\}$ is biorthogonal to $\{\widetilde{\phi}_\alpha (t-n) : n \in \mathbb Z\}$ if and only if 
$$\sum_{k\in \mathbb Z} \Theta_\alpha(u +2k\pi\sin\alpha)\,\overline{\widetilde{\Theta}_\alpha(u+2k\pi\sin\alpha)} = \dfrac{1}{\sin\alpha}$$
where $\Theta_\alpha(u)$ is FrFT of $\phi(t)$.

\parindent=0mm \vspace{.1in}
\textbf{Proof.}
\begin{align*}
\mathcal{F}_\alpha \left(\phi_{\alpha,n}(t)\right) (u)& =\int_{-\infty}^{\infty} \phi(t-n)e^{-j/2(t^2-n^2-(t-n)^2)\cot\alpha}\mathcal{K_\alpha}(u,t)\,du\\\
& = \mathcal{A_\alpha}\int_{-\infty}^{\infty} \phi(t-n) e^{j/2 [(t-n)^2+ u^2]\cot\alpha - j t u\csc\alpha+jn^2/2\cot\alpha}\,du \\\
& = e^{jn^2/2\cot\alpha-jnu\csc\alpha}\mathcal{A_\alpha}\int_{-\infty}^{\infty} \psi(t-n)e^{j/2[(t-n)^2+u^2]\cot\alpha-j(t-n)u\csc\alpha}\,du \\\
& =e^{jn^2/2 \cot\alpha-jnu\csc\alpha}\Theta_\alpha(u).\tag{3.1}
\end{align*}
Now from the Parseval identity of FrFT, we have
\begin{align*}
\left\langle\phi_{\alpha,n}(t) \widetilde{\phi}_{\alpha,m}(t) \right\rangle & =\left\langle\mathcal{F_\alpha}(\phi_{\alpha,n}(t))(u),\mathcal{F_\alpha}(\phi_{\alpha,m}(t))(u) \right\rangle \\\
&= \left \langle e^{\frac{jn^2}{2} \cot\alpha-jnu\csc\alpha}\Theta_\alpha(u),e^{\frac{jm^2}{2} \cot\alpha-jmu\csc\alpha}\widetilde{\Theta}_\alpha(u) \right\rangle \\\
& =\int_{-\infty}^{\infty}e^{\frac{j}{2}(n^2-m^2)\cot\alpha - j(n-m)u\csc\alpha}\Theta_\alpha(u)\overline{\widetilde{\Theta_\alpha}}(u)\,du.\tag{3.2}
\end{align*}

\parindent=0mm \vspace{.1in}
Since
$$\left\langle \phi_{\alpha,n}(t), \widetilde{\phi}_{\alpha,m}(t) \right \rangle = \delta_{n,m} ~~\forall~~ n,m\eqno(3.3)$$
Setting $n = n-m,$ then from $(3.2)$ and $(3.3)$ it follows that 
$$\int_{-\infty}^{\infty} e^{-jnu\csc\alpha}\Theta_\alpha(u)\overline{\widetilde{\Theta}_\alpha (u)} =\delta_{n,0} \eqno(3.4)$$
implies
\begin{align*}
\sum_{k\in \mathbb Z}&\int_{0}^{2\pi\sin\alpha} e^{-jn(u+2k\pi\sin\alpha)\csc\alpha}\Theta_\alpha(u+2k\pi\sin\alpha)\overline{\widetilde{\Theta}_\alpha (u+2k\pi\sin\alpha)}\, du~~~~~~~~~~~~~~~~~~~~~~~ \\\
&\qquad\qquad\qquad =\int_{0}^{2\pi\sin\alpha} e^{-jnu\csc\alpha} \sum_{k \in \mathbb Z}\Theta_\alpha(u+2k\pi\sin\alpha)\overline{\widetilde{\Theta}_\alpha (u+2k\pi\sin\alpha)}\, du\\\
&\qquad\qquad\qquad = \delta_{n,0}.\tag{3.5}
\end{align*}

Let
$$ \mathcal{L}(u) = \sum_{k \in \mathbb Z}\Theta_\alpha(u+2k\pi\sin\alpha)\overline{\widetilde{\Theta}_\alpha (u +2k\pi\sin\alpha)},$$
then we have
$$ \mathcal{L}(u+2\pi\sin\alpha) = \sum_{k \in \mathbb Z}\Theta_\alpha(u+2(k+1)\pi\sin\alpha)\overline{\widetilde{\Theta}_\alpha (u +2(k+1)\pi\sin\alpha)}.\eqno(2.6)$$
By setting $k' = k+1$ in $(3.6)$, we obtain that 
\begin{align*}
\mathcal{L}(u)& = \sum_{k' \in \mathbb Z}\Theta_\alpha(u+2k'\pi\sin\alpha)\overline{\widetilde{\Theta}_\alpha (u +2k'\pi\sin\alpha)} \\\
&= \mathcal{L}(u).\tag{3.7}
\end{align*}
It clearly implies that $\mathcal{L}(u)$ is $2k\pi\sin\alpha$ periodic function. Therefore from $(3.7)$, we have
$$\int_{0}^{2\pi\sin\alpha} e^{-jnu\csc\alpha}\mathcal{L}(u)\, du = \delta_{n,0}  = \delta_{n\csc\alpha,0}. \eqno(3.8)$$

If we set $n' = n \csc\alpha$ in $(3.8)$, we obtain
$$\int_{0}^{2\pi\sin\alpha}\mathcal{L}(u)e^{-jn'u}\,du= \delta_{n',0},\eqno(3.9)$$
implies,
$$\dfrac{1}{2\pi \sin\alpha}\int_{0}^{2\pi\sin\alpha}\mathcal{L}(u)e^{-jn'u}\,du= \dfrac{1}{2\pi \sin\alpha}\delta_{n',0},$$
which further implies,
$$ \mathcal{L}(u) = \mathcal{F}^{-1}\left\{\dfrac{1}{2\pi \sin\alpha}\delta_{n',0}\right\}(u) = \dfrac{1}{\sin\alpha}.$$
Thus we have
$$\sum_{k \in \mathbb Z}\Theta_\alpha(u+2k\pi\sin\alpha)\overline{\widetilde{\Theta}_\alpha (u +2k\pi\sin\alpha)}= \dfrac{1}{\sin\alpha}.\eqno\square$$

The following result etablishes a sufficient condition for the translates of a function to be linearly independent.

\parindent=0mm \vspace{.1in} 
\textbf{Lemma 3.6.} Let $\phi_\alpha(t) \in L^2(\mathbb R)$. Assume that there exist constants $ C_1, C_2 >0$ such that

$$C_1 \le \sum_{k \in \mathbb N} \left|\Theta_\alpha(u+2k\pi\sin\alpha)\right|^2 \le C_2~~\forall ~~u \in \mathbb R. \eqno(3.10)$$
Then $\{\phi_\alpha(t-n) : n \in \mathbb Z\}$ is linearly independent.

\parindent=0mm \vspace{.1in} 
\textbf{Proof.} By virtue of Lemma 3.3, it suffices to find a function $\widetilde{\phi}_\alpha(t)$ whose translates are biorthogonal to the translates of $\phi_\alpha.$ We define $\widetilde{\phi}_\alpha$ by
$$ \mathcal{F}_\alpha\{\widetilde{\phi}_\alpha(t)\} (u) = \dfrac{\Theta_\alpha(u)}{\displaystyle\sum_{k \in \mathbb Z}\left|\Theta_\alpha(u+2k\pi\sin\alpha)\right|^2}.$$

By virtue of $(3.10)$, this function is well defined. Now
\begin{align*}
\sum_{m \in \mathbb N}\Theta_\alpha(u+2m\pi\sin\alpha)&\overline{\widetilde{\Theta}_\alpha(u+2k\pi\sin\alpha)}~~~~~~~~~~~~~~~~~~~~~~\\\
&=\sum_{m \in \mathbb N}\Theta_\alpha(u+2m\pi\sin\alpha)\dfrac{\overline{\Theta_\alpha(u+2k\pi\sin\alpha)}}{\displaystyle\sum_{k \in \mathbb N}\left|\Theta_\alpha(u+2k\pi\sin\alpha+2m\pi\sin\alpha)\right|^2} \\\
&=\dfrac{\displaystyle\sum_{m \in \mathbb N}\left|\Theta_\alpha(u+2m\pi\sin\alpha)\right|^2}{\displaystyle\sum_{\ell \in \mathbb N}\left|\Theta_\alpha(u+2\ell\pi\sin\alpha)\right|^2} \\\
& = 1.
\end{align*}
It clearly implies $\{\phi_\alpha(t-n)\}$ is biorthogonal to $\{\widetilde{\phi}_\alpha(t-n)\}.~~~~\square.$

\parindent=0mm \vspace{.1in}
\textbf{Lemmma 3.7.} Suppose that $\phi_\alpha(t)$ satisfies $(3 .10)$. Any function $f$ in span $\{\phi_\alpha(t-n) : n \in \mathbb Z\}$ is of the form
$$ f = \sum_{n \in \mathbb N}a[n]\phi_\alpha(t-n),$$
where $a[n] \in \ell^2(\mathbb Z)$ is a finite sequence. Let $\widetilde{a}_\alpha(u)$ be the discrete FrFT of $a[n]$. Then
$$C_1\int_{0}^{2\pi\sin\alpha} |\widetilde{a}_\alpha(u)|^2\,du \le \|f\|_2^2 \le C_2 \int_{0}^{2\pi \sin \alpha} |\widetilde{a}_\alpha(u)|^2\,du.$$

\parindent=0mm \vspace{.1in}
\textbf{Proof.} By Plancherel theorem, we have
\begin{align*}
\int_{-\infty}^{\infty}|f(t)|^2\,dt & = \int_{-\infty}^{\infty}\left|\sum_{n\in \mathbb N} a[n] \phi_\alpha(t-n)\right|^2\,dt \\\
&=\int_{-\infty}^{\infty}\left|\Theta_\alpha(u)\right|^2 \,|\widetilde{a}_\alpha[u]|^2\,du \\\
&=\int_{0}^{2\pi \sin\alpha} \sum_{k \in \mathbb N} \left|\Theta_\alpha(u+2k\pi\sin\alpha)\widetilde{a}_\alpha[u]\right|^2\,du.
\end{align*}
By invoking $(3.10)$, the result follows.

\parindent=0mm \vspace{.1in}

\textbf{Lemma 3.8.} Let $\{\phi_\alpha(t-n) : n \in \mathbb Z\}$ be a Reisz basis for its closed linear span. Suppose that there exists a function $\widetilde{\phi}_\alpha$ such that $\{\widetilde{\phi}_\alpha(t-n) :n \in \mathbb Z\}$ is biorthogonal to $\{\phi_\alpha(t-n) : n\in \mathbb Z\}$. Then

\parindent=0mm \vspace{.1in}
(a) for every $f \in \overline{span}\{\phi_\alpha(t-n) : n\in \mathbb Z\},$ we have
$$ f = \sum_{n \in \mathbb Z}\left\langle f,\widetilde{\phi}_\alpha(t-n)\right\rangle \phi_\alpha(t-n);\eqno(3.11)$$

\parindent=0mm \vspace{.1in}
(b) there exist constants $ A,B >0$ such that for every $ f \in \overline{span}\{\phi_\alpha(t-n) : n \in \mathbb Z\},$ we have
$$A\|f\|_2^2 \le \sum_{n=1}^{\infty}\left|\left\langle f,\widetilde{\phi}_\alpha(t-n)\right\rangle\right|^2 \le B\|f\|_2^2.\eqno(3.12)$$

\parindent=0mm \vspace{.1in}
\textbf{Proof.} Since $\left\{\phi_\alpha(t-n) : n \in \mathbb Z\right\}$ forms a Riesz basis for its closed linear span, then there exist constants $C_1$ and $C_2$ such that $(3.10)$ holds. First establish the results for $f \in span\{\phi_\alpha(t-n) : n\in \mathbb Z\}$ and we generalize the established results to $ \overline{span}\left\{\phi_\alpha(t-n) : n\in \mathbb Z\right\}$.

\parindent=8mm \vspace{.1in}
Let $f \in span\left\{\phi_\alpha(t-n) : n\in \mathbb Z\right\}$, then there exists a finite sequence $a[n]$ such that
$$f = \sum_{n=1}^{\infty} a[n] \phi_\alpha(t-n).$$
By the definition of biorthogonality, we have
\begin{align*}
\left\langle f, \widetilde{\phi}_\alpha(t-k)\right\rangle &= \left\langle \sum_{n=1}^{\infty} a[n] \phi_\alpha(t-n), \widetilde{\phi}_\alpha(t-k)\right\rangle\\\
&=\sum_{n=1}^{\infty} a[n]\left\langle  \phi_\alpha(t-n), \widetilde{\phi}_\alpha(t-k)\right\rangle\\\
& = a[k].
\end{align*}
Thus $(a)$ is established for $f \in span\left\{\phi_\alpha(t-n) : n\in \mathbb Z\right\}$.

\parindent=8mm \vspace{.1in}
Now we proceed to establish $(b)$. Since $(3.10)$ is satisfied, by Lemma 2.7, for every $f \in span\left\{\phi_\alpha(t-n) : n\in \mathbb Z\right\}$, we have
$$C_2^{-1}\|f\|_2^2 \le \int_{0}^{2\pi\sin\alpha} \left|\widetilde{a}_\alpha[u]\right|^2\,du \le C_1^{-1} \|f\|_2^2.$$
By Plancherel formula for the Fourier series and the fact $a[n] = \langle f, \widetilde{\phi}_\alpha(t-n)\rangle,$ we have
$$ \int_{0}^{2\pi\sin \alpha}\left|\widetilde{a}_\alpha[u]\right|^2\,du = \sum_{n \in \mathbb N}|a[n]|^2 = \sum_{n \in \mathbb N} \left|\left\langle f, \widetilde{\phi}_\alpha(t-n)\right\rangle\right|^2,$$
Thus $(b)$ is obtained.

\parindent=8mm \vspace{.1in}
Finally we proceed to generalize the results to $\overline{span}\{\phi_\alpha(t-n) : n\in \mathbb Z\}$. we first establish $(b)$. For $ f \in \overline{span}\{\phi_\alpha(t-n) : n\in \mathbb Z\}$, there exists a sequence $f[m] \in \ell^2(\mathbb N)$ in $span\{\phi_\alpha(t-n) : n\in \mathbb Z\}$ such that $\lim_{m\rightarrow \infty} f[m] = f.$ Hence for each $n \in \mathbb N$, we have
$$\left\langle f[m], \widetilde{\phi}_\alpha(t-n)\right \rangle \rightarrow \left\langle f ,\widetilde{\phi}_\alpha(t-n)\right\rangle~~\textit{as}~~ m \rightarrow \infty.$$
The result holds for each $f[m]$. Hence,
\begin{align*}
\sum_{n=-N}^{N} \left|\left\langle f , \widetilde{\phi}_\alpha(t-n)\right\rangle\right|^2 &=\sum_{n=-N}^{N} \lim_{m \rightarrow \infty}\left|\left\langle f[m] , \widetilde{\phi}_\alpha(t-n)\right\rangle\right|^2\\\
&=\lim_{m \rightarrow \infty}\sum_{n=-N}^{N} \left|\left\langle f[m] , \widetilde{\phi}_\alpha(t-n)\right\rangle\right|^2\\\
&\le B\lim_{m \rightarrow \infty}\|f[m]\|_2^2\\\
&=B \|f\|_2^2.
\end{align*} 
Letting $N \rightarrow \infty$ in the above expression, we obtain
$$\sum_{n\in \mathbb Z} \left|\left\langle f , \widetilde{\phi}_\alpha(t-n)\right\rangle\right|^2  \le B \|f\|_2^2$$
Hence the upper bound in $(3.12)$ holds. Now  by the Cauchy Schwarz inequality for sequences, we have for each $m \in \mathbb Z$
$$\left\{\sum_{n \in \mathbb Z} \left|\left\langle f[m] , \widetilde{\phi}_\alpha(t-n)\right\rangle\right|^2\right\}^{1/2} \le \left\{\sum_{n \in \mathbb Z} \left|\left\langle f[m]-f , \widetilde{\phi}_\alpha(t-n)\right\rangle\right|^2\right\}^{1/2}+\left\{\sum_{n \in \mathbb Z} \left|\left\langle f , \widetilde{\phi}_\alpha(t-n)\right\rangle\right|^2\right\}^{1/2}$$

\parindent=0mm \vspace{.1in}
Since the upper bound in $(3.12)$ holds for each $f[m] - f$ and the lower bound holds for each $f[m]$, we have
$$ A^{1/2} \|f[m]\|_2 \le B^{1/2} \|f[m] -f\|_2 +\left\{\sum_{n \in \mathbb Z} \left|\left\langle f , \widetilde{\phi}_\alpha(t-n)\right\rangle\right|^2\right\}^{1/2}$$
Taking limits as $m \rightarrow \infty,$ we get
$$A\|f\|_2^2 \le \sum_{n \in \mathbb N} \left|\left\langle f, \widetilde{\phi}_\alpha(t-n)\right\rangle\right|^2$$
which is the upper bound in $(3.12)$.

\parindent=8mm \vspace{.1in}
Now we will prove $(a)$ for $f \in \overline{span}\{\phi_\alpha(t-n) : n\in \mathbb Z\}$. Let $\epsilon >0$ and $g \in span\{\phi_\alpha(t-n) : n\in \mathbb Z\}$ such that $\|f-g\|_2 < \epsilon.$ Since $(a)$ holds for every $g$, therefore for large $N, M \in \mathbb N,$ we have
\begin{align*}
f - &\sum_{n=-M}^{N} \left\langle f , \widetilde{\phi}_\alpha(t-n)\right\rangle \phi_\alpha(t-n)~~~~~~~~~~\\\
&= f -g +\sum_{n=-M}^{N} \left\langle g , \widetilde{\phi}_\alpha(t-n)\right\rangle \phi_\alpha(t-n)-\sum_{n=-M}^{N} \left\langle f , \widetilde{\phi}_\alpha(t-n)\right\rangle \phi_\alpha(t-n) \\\
&=f-g +\sum_{n=-M}^{N} \left\langle g- f , \widetilde{\phi}_\alpha(t-n)\right\rangle \phi_\alpha(t-n).
\end{align*}
Hence, by Cauchy Schwarz Inequality, we have
\begin{align*}
\left\|f - \sum_{n=-M}^{N} \left\langle f , \widetilde{\phi}_\alpha(t-n)\right\rangle \phi_\alpha(t-n)\right\|_2 &\le \|f-g\|_2 + \left\|\sum_{n=-M}^{N} \left\langle g- f , \widetilde{\phi}_\alpha(t-n)\right\rangle \phi_\alpha(t-n)\right\| \\\
&\le \|f-g\|_2 +\sqrt{C_2}\left\{\sum_{n=-M}^{N}\left|\left\langle g- f , \widetilde{\phi}_\alpha(t-n)\right\rangle\right|^2\right\}^{1/2} \\\
&\le \|f-g\|_2 +\sqrt{C_2}\sqrt{B}\|f-g\|_2 \\\
& <(1+\sqrt{C_2B})\epsilon.
\end{align*}
Since $\epsilon $ is arbitrary, the result follows.$\square$

\parindent=0mm \vspace{.1in}
\textbf{4. Fractional Multiresolution Analysis Associated with Fractional Wavelets}

\parindent=0mm \vspace{.1in}
As in case of conventional wavelets, there corresponds a multiresolution analysis. In the similar way fractional wavelets give rise to fractional multiresolution analysis.

\parindent=0mm \vspace{.1in}
\textbf{Definition 4.1.} An MRA associated with the fractional wavelet transform is defined as a sequence of closed subspaces $\{V_k^\alpha\} \in L^2(\mathbb R)$ such that

\parindent=0mm \vspace{.1in}
(i) $V_k^\alpha \subseteq V_{k+1}^\alpha,~~k \in \mathbb Z;$

\parindent=0mm \vspace{.1in}
(ii) $\bigcup_{k \in \mathbb Z} V_k^\alpha$ is dense in $L^2(\mathbb R)$;

\parindent=0mm \vspace{.1in}
(iii) $\bigcap_{k \in \mathbb Z} V_k^\alpha = \{0\}$;

\parindent=0mm \vspace{.1in}
(iv) $f(t) \in V_k^\alpha$ if and only if $f(2t)e^{\frac{j}{2}[(2t)^2 - t^2]\cot\alpha} \in V_{k+1}^\alpha,~~k \in \mathbb Z;$

\parindent=0mm \vspace{.1in}
(v) there is a function $ \phi \in V_0^\alpha $ called {\it scaling function} such that $\{\phi_{\alpha,0,n} = \phi(t-n)e^{-j(tn+n^2)\cot\alpha} : n \in \mathbb Z\}$ is an orthonormal basis of subspace $V_0^\alpha$.

\parindent=8mm \vspace{.1in}
In the above definition, if we assume that the set of functions $\{\phi_{\alpha,0,n} : n\in \mathbb Z\}$ form a Reisz basis of $V_0^\alpha$, then $\phi(t)$ generates a generalized fractional MRA $\{V_m^\alpha\}$ of $L^2(\mathbb R)$, then
$$\phi_{\alpha, m,n}(t) = 2^{\frac{m}{2}}\phi(2^m t-n)e^{\frac{-j}{2}[t^2-(2^{-m}n)^2 -(2^mt-n)^2]\cot\alpha}$$
is the orthonormal basis of $\{V_m^\alpha\}$.

\parindent=0mm \vspace{.1in}
\textbf{Theorem 4.2.} Let $\phi \in L^2(\mathbb R)$ such that the collection  $\{\phi_{\alpha,0,n} (t) : n \in \mathbb Z\}$ is a Reisz basis of the space
$$V_0^\alpha = \left\{\sum_{n \in\mathbb Z} c[n] \phi_{\alpha,0,n}(t) : c[n] \in \ell^2(\mathbb Z)\right\}$$
of $L^2(\mathbb R)$ if and only if there exists positive constants $ A, B$ such that for all $ u\in I = [0, 2\pi\sin\alpha]$, we have

\parindent=0mm \vspace{.1in}
$$ A \le \mathcal{G}^2(\alpha, \phi,u) \le B \eqno(4.1)$$
where

\parindent=0mm \vspace{.1in}
$$\mathcal{G}(\alpha, \phi ,u) = \sqrt{2\pi\sin\alpha\sum_{k \in \mathbb Z}\left|\Theta_\alpha(u+2k\pi\sin\alpha)\right|^2}.\eqno(4.2)$$

\parindent=0mm \vspace{.1in}
\textbf{Proof.} For any $f(t) \in V_0^\alpha$, we have
$$ f(t) = \sum_{n\in \mathbb Z} c[n] \phi_{\alpha,0,n}(t)\eqno(4.3)$$
where $c[n]\in\ell^2(\mathbb Z)$.

\parindent=0mm \vspace{.1in}
On taking FrFT on both sides of $(4.3)$, we obtain
$$ \mathcal{F}_\alpha\{f(t)\}(u) = \sqrt{2\pi}\, \widetilde{c}_\alpha(u)\Theta_\alpha(u)\eqno(4.4)$$
where $\widetilde{c}_\alpha(u)$ denotes the discrete FrFT of $c[n]$. By using Parseval formula of the FrFT, we have
\begin{align*}
\|f(t)\|_{L^2(\mathbb R)}^2 & = \|\mathcal{F}_\alpha\{f(t)\}(u)\|_{L^2(\mathbb R)}^2\\\
&= \int_{-\infty}^{\infty} 2\pi |\widetilde{c}_\alpha(u)|^2\,|\Theta_\alpha(u)|^2\,du\\\
&= \sum_{k \in \mathbb Z}\int_{0}^{2\pi\sin\alpha}2\pi |\widetilde{c}_\alpha(u+2k\pi\sin\alpha)|^2\,|\Theta_\alpha(u+2k\pi\sin\alpha)|^2\,du\\\
& =\int_{0}^{2\pi\sin\alpha}|\widetilde{c}_\alpha(u)|^2\mathcal{G}^2(\alpha, \phi,u)\,du,\tag{4.5}
\end{align*}
Further, Parsevals formula for discrete FrFT yields
$$\|c[n]\|_{\ell^2(\mathbb Z)}^2 = \sum_{n \in \mathbb Z}|c[n]|^2 = \int_{0}^{2\pi\sin\alpha}|\widetilde{c}_\alpha(u)|^2\,du.\eqno(4.6)$$

Now, Eqns. $(4.1), (4.2) $ and $(4.3)$ yields
$$ A \|c[n]\|_{\ell^2(\mathbb Z)}^2 \le \left\|\sum_{n \in \mathbb Z}c[n]\phi_{\alpha,0,n}(t)\right\|^2 \le B \|c[n]\|_{\ell^2(\mathbb Z)}^2.\eqno(4.7)$$

It follows from $(4.7)$ and Definition 3.4. that $\{\phi_{\alpha,0,n}(t) : n \in \mathbb Z\}$ is a Reisz basis for $V_0^\alpha$. In particular $\{\phi_{\alpha,0,n}(t) : n \in \mathbb Z\}$ is an orthonormal basis for $V_0^\alpha$ if and only if $ A = B = 1.$

\parindent=0mm \vspace{.1in}
\textbf{Lemma 4.3.} Let $\phi$ be the scaling function for fractional MRA $\{V_j^\alpha : j \in \mathbb Z\}.$ Then for each $ j \in \mathbb Z,~\{\phi_{\alpha,j,k} : k \in \mathbb Z\}$ is a Reisz basis for $V_j^\alpha$.

\parindent=0mm \vspace{.1in}
\textbf{Proof.} Define $\widetilde{\phi}$ by 

$$ \mathcal{F}_\alpha\{\widetilde{\phi}(t)\}(u) = \dfrac{\Theta_\alpha (u)}{\sqrt{2\pi\sin\alpha\sum_{k \in \mathbb Z}\left|\Theta_\alpha(u+2k\pi\sin\alpha)\right|^2}},$$
then by the same arguments as in the proof of Lemma 3.4, $\{\widetilde{\phi}_\alpha(t- k) : k \in \mathbb Z\}$ is biorthogonal to $\{\phi_\alpha(t- k) : k \in \mathbb Z\}$

\parindent=0mm \vspace{.1in}
Hence,
\begin{align*}
\left\langle \phi_{\alpha,j,n}, \widetilde{\phi}_{\alpha, j,m}\right\rangle &= \left\langle \delta_j\phi_\alpha(t-n), \delta_j \widetilde{\phi}_\alpha(t-m)\right\rangle\\\
& = \left\langle \phi_\alpha(t-n), \widetilde{\phi}_\alpha(t-m)\right\rangle = \delta _{n,m}
\end{align*}
that is, $\left\{\widetilde{\phi}_{\alpha, j,k} : k \in \mathbb Z\right\}$ is biorthogonal to $\left\{\phi_{\alpha , j,k} : k \in \mathbb Z\right\}$ for every $j \in \mathbb Z.$ Therefore by Lemma 3.3. $\left\{\phi_{\alpha , j ,k} : k \in \mathbb Z\right\}$ is linearly independent.

\parindent=8mm \vspace{.1in}
Now we need to show  the collection $\{\phi_{\alpha, j,k} : k \in \mathbb Z\}$ satisfies the frame condition. For any $f \in V_j^\alpha,$ we have

\begin{align*}
\sum_{k \in \mathbb Z} \left| \left \langle f , \phi_{\alpha, j,k}\right\rangle \right|^2 &= \sum_{k \in \mathbb Z} \left| \left \langle f , \delta_j\phi_{\alpha}(t-k)\right\rangle \right|^2\\\
& = \sum_{k \in \mathbb Z} \left| \left \langle\delta_{-j} f ,\phi_{\alpha}(t-k)\right\rangle \right|^2.
\end{align*}
Since $\{\phi_\alpha (t-n) : k \in \mathbb Z\}$ is a Riesz basis for $V_0^\alpha$ and $\delta_{-j}f \in V_0^\alpha,$ there exist constants $ A, B > 0$ such that for every $f \in V_j^\alpha$,

$$ A\|\delta_{-j}f\|_2^2 \le \sum_{k \in \mathbb Z} \left|\left\langle \delta_{-j} f, \phi_\alpha(t-k)\right\rangle\right|^2 \le B\|\delta_{-j}f\|_2^2.$$
This is equivalent to 

$$ A\|f\|_2^2 \le \sum_{k \in \mathbb Z} \left|\left\langle f, \phi_\alpha(t-k)\right\rangle\right|^2 \le B\|f\|_2^2.$$
Hence, $\{\phi_{\alpha , j,k} : k \in \mathbb Z\}$ satisfies the frame condition.

\parindent=0mm \vspace{.1in}
\textbf{Lemma 4.4.} Suppose that $\left\{ V_j^\alpha : j\in \mathbb Z\right\}$ is a fractional MRA with scaling function $\phi$. Then there exists a sequence $\{h[n] \}_{n \in \mathbb Z}$ in $\ell^2 (\mathbb Z)$ called the scaling filter such that

$$\phi(t) = \sum_{n \in \mathbb Z} h[n] \sqrt{2}\phi(2t-n)e^{\frac{-j}{2}[t^2-(\frac{n}{2})^2 - (2t-n)^2]\cot\alpha}$$
and a $2k\pi\sin\alpha$-periodic function $\Lambda_\alpha$  called the $\it {auxilliary~ function}$ such that
$$\Theta_\alpha(u) = \Lambda_\alpha\left(\frac{u}{2}\right)\Theta_\alpha\left(\frac{u}{2}\right).$$

\textbf{Proof.} By Lemma  4.3, $\left\{\phi_{\alpha ,1,n} :n \in \mathbb Z\right\}$ is a Riesz basis for $V_1^\alpha$. Since $ \phi_{\alpha,0,0} (t) \in V_0^\alpha \subseteq V_1^\alpha,$ therefore by virtue of $(3.11)$ there must exist a sequence $\{h[n]\}_{n \in \mathbb Z}$ in $\ell^2(\mathbb Z)$ such that

$$ \phi_{\alpha,0,0} = \sum_{n \in \mathbb Z} h[n] \phi_{\alpha,1,n}(t)$$
which can be simplified as
$$ \phi(t) = \sum_{n \in \mathbb Z} h[n] \sqrt{2}\phi(2t-n)e^{\frac{-j}{2}[t^2-(\frac{n}{2})^2 - (2t-n)^2]\cot\alpha} \eqno(4.8)$$
and the coefficient can be solved as 
$$ h[n] = \sqrt{2}\int_{-\infty}^{\infty} \phi(t) \phi^{*}(2t-n)e^{\frac{j}{2}[t^2-(\frac{n}{2})^2 - (2t-n)^2]\cot\alpha}dt. \eqno(4.9)$$
By taking the FrFT on both sides of Eq. (4.8), we have
\begin{align*}
\Theta_\alpha(u) & =\sum_{n \in \mathbb Z} h[n]\sqrt{2}\mathcal{A}_\alpha\int_{-\infty}^{\infty}\phi(2t-n)e^{\frac{j}{2}[t^2-(\frac{n}{2})^2 - (2t-n)^2]\cot\alpha -jtu\csc\alpha}dt \\\
&=\dfrac{1}{\sqrt{2}}e^{\frac{3ju^2}{8}\cot\alpha}\sum_{n \in \mathbb Z} h[n]e^{\frac{jn^2}{8}\cot\alpha -\frac{jnu}{2}\csc\alpha}\\\
&\qquad\qquad\qquad\times\mathcal{A}_\alpha\int_{-\infty}^{\infty}\phi(2t-n)e^{\frac{j}{2}[t^2-(\frac{n}{2})^2 - (2t-n)^2]\cot\alpha -j(2t-n)u\csc\alpha}d(2t-n) \\\
&=\dfrac{1}{\sqrt{2}}e^{\frac{3ju^2}{8}\cot\alpha}\sum_{n \in \mathbb Z} h[n]e^{\frac{jn^2}{8}\cot\alpha -\frac{jnu}{2}\csc\alpha}\Theta_\alpha\left(\frac{u}{2}\right)\\\
&=\dfrac{1}{\sqrt{2}}e^{\frac{3ju^2}{8}\cot\alpha}D_\alpha\left(\frac{u}{2}\right)\Theta_\alpha\left(\frac{u}{2}\right),\tag{4.10}
\end{align*}
where
$$D_\alpha(u) =\sum_{n \in \mathbb Z} h[n]e^{\frac{jn^2}{8}\cot\alpha -\frac{jnu}{2}\csc\alpha}.$$
By defining
\begin{align*}
\Lambda_\alpha(u) &=\dfrac{1}{\sqrt{2}}e^{\frac{3ju^2}{2}\cot\alpha}D_\alpha(u) \\\
& =\dfrac{1}{\sqrt{2}}\sum_{n \in \mathbb Z} f[n] \mathcal{A}_\alpha e^{\frac{jn^2}{2}\cot\alpha - jnu\csc\alpha}.
\end{align*}
Eq.(4.10) can be written as 
$$\Theta_\alpha(u) = \Lambda_\alpha\left(\frac{u}{2}\right)\Theta_\alpha\left(\frac{u}{2}\right).$$

It is to be noted that $\Lambda_\alpha(u) $ is  a $2k\pi\sin\alpha$-periodic function since we have
\begin{align*}
\Lambda_\alpha(u+2k\pi\sin\alpha) &= \dfrac{1}{\sqrt{2}}\sum_{n \in \mathbb Z} f[n] \mathcal{A}_\alpha e^{\frac{jn^2}{2}\cot\alpha - jn(u +2k\pi\sin\alpha)\csc\alpha} \\\
&=\dfrac{1}{\sqrt{2}}\sum_{n \in \mathbb Z} f[n] \mathcal{A}_\alpha e^{\frac{jn^2}{2}\cot\alpha - jnu\csc\alpha}\\\
&= \Lambda_\alpha(u).
\end{align*}

\parindent=0mm \vspace{.1in}
\textbf{Definition 4.5.} A pair of fractional MRAs $\left\{V_j^\alpha : j\in \mathbb Z\right\}$ and $\left\{\widetilde{V}_j^\alpha : j\in \mathbb Z\right\}$ with scaling functions $\phi$ and $\widetilde{\phi}$ respectively are said to be dual to each other if $\left\{\phi_\alpha(t-k) : k\in \mathbb Z\right\}$ and $\left\{\widetilde{\phi}_\alpha(t-k) : k\in \mathbb Z\right\}$ are biorthogonal.

\parindent=0mm \vspace{.1in}
\textbf{Definition 4.6.} Let $\phi$ and $\widetilde{\phi}$ be scaling functions for dual MRAs. For each $j \in \mathbb Z$ we define the operators $\mathcal{P}_{\alpha,j},~\widetilde{\mathcal{P}}_{\alpha,j}$ on $L^2(\mathbb R)$ by 
$$\mathcal{P}_{\alpha,j}f = \sum_{k\in \mathbb Z}\left\langle f , \widetilde{\phi}_{\alpha,j,k}\right\rangle \phi_{\alpha,j,k},$$

$$\widetilde{\mathcal{P}}_{\alpha,j} =\sum_{k\in \mathbb Z}\left\langle f , \phi_{\alpha,j,k}\right\rangle \widetilde{\phi}_{\alpha,j,k}.$$

\textbf{Lemma 4.7.} The operators $\mathcal{P}_{\alpha,j},~\widetilde{\mathcal{P}}_{\alpha,j}$ are uniformly bounded.

\parindent=0mm \vspace{.1in}
\textbf{Proof.} Since the translates of $\phi$ and $\widetilde{\phi}$ form Riesz basis for their closed linear spans, therefor by Lemma 4.2, there exist constants $C_1$ and $C_2$ such that

\parindent=0mm \vspace{.0in}
$$ c_1 \le \sum_{k\in\mathbb Z} \left|\Theta_\alpha(u+2k\pi\sin\alpha)\right|^2 \le C_2$$
and
$$ c_1 \le \sum_{k\in\mathbb Z} \left|\widetilde{\Theta}_\alpha(u+2k\pi\sin\alpha)\right|^2 \le C_2.$$

For a sequence $\{c[k]\}_{k \in \mathbb Z} \in \ell^2(\mathbb Z)$, there exists $B >0$ such that

\parindent=0mm \vspace{.0in}
$$ \left\|\sum_{k \in \mathbb Z}c[k]\phi_{\alpha,0,k}\right\|_2^2 \le \sum_{k \in \mathbb Z}|c[k]|^2.$$

\parindent=0mm \vspace{.0in}
Now for any $f \in L^2(\mathbb R)$, we have 
\begin{align*}
\sum_{k \in \mathbb Z} \left|\left\langle f, \phi_{\alpha,0,k}\right\rangle \right|^2 &= \sum_{k\in \mathbb Z} \left|\int_{-\infty}^{\infty}\mathcal{F}_\alpha\{f(t)\}(u) \overline{\Theta_\alpha(u)}e^{\frac{jn^2}{2}\cot\alpha-jnu\csc\alpha}\,du\right|^2\\\
& = \int_{0}^{2\pi\sin\alpha}\left|\sum_{\ell\in \mathbb Z}\mathcal{F}_\alpha\{f(t)\}(u+2\pi\ell\sin\alpha) \overline{\Theta}_\alpha(u+2\pi\ell\sin\alpha)\,du\right|^2\\ \
&  \le \int_{0}^{2\pi\sin\alpha}\left\{\sum_{\ell\in \mathbb Z}\left|\mathcal{F}_\alpha\{f(t)\}(u+2\pi\ell\sin\alpha)\right|^2\right\} \left\{\sum_{\ell\in \mathbb Z}\left|\Theta_\alpha(u+2\pi\ell\sin\alpha)\right|^2\right\}\,du\\\
& \le \int_{0}^{2\pi\sin\alpha}\sum_{\ell\in \mathbb Z}\left|\mathcal{F}_\alpha\{f(t)\}(u+2\pi\ell\sin\alpha)\right|^2\,du \\\
& =C_2 \int_{-\infty}^{\infty} \left|\mathcal{F}_\alpha\{f(t)\}(u)\right|^2\,du\\\
& = C_2\|f\|_2^2.
\end{align*}

\parindent=0mm \vspace{.0in}
Similar estimates hold for $\widetilde{\phi}$. Hence for $f \in L^2(\mathbb R)$, we  have

\parindent=0mm \vspace{.0in}
\begin{align*}
\|\mathcal{P}_{\alpha,0} f \|_2^2 &= \left\|\sum_{k \in \mathbb Z} \left\langle f, \widetilde{\phi}_{\alpha,0,k}\right\rangle \phi_{\alpha,0,k}\right\|_2^2 \\ \
& \le B \sum_{k \in \mathbb Z}\left|\left\langle f,\widetilde{\phi}_{\alpha,0,k}\right\rangle\right|^2\\\
& \le B C_2\|f\|_2^2.
\end{align*}

\parindent=0mm \vspace{.0in}
Thus $\mathcal{P}_{\alpha,0}$ is a bounded operator on $L^2(\mathbb R)$ with norm at most $\sqrt{BC_2} =C,$(say). Since the dilation operators are unitary and since
\begin{align*}
\mathcal{P}_{\alpha,j} f &= \sum_{k \in \mathbb Z} \left\langle f, \widetilde{\phi}_{\alpha,j,k}\right\rangle \phi_{\alpha,j,k}\\\
& =\sum_{k \in \mathbb Z} \left\langle \delta_{-j}f, \widetilde{\phi}_{\alpha,0,k}\right\rangle \delta_{-j}\phi_{\alpha,0,k},
\end{align*}
we conclude that the operator norm of $\mathcal{P}_{\alpha,j}$ is at most $C$. Similar arguments work for $\widetilde{\mathcal{P}}_{\alpha,j}$. This finishes the proof of the lemma.

\parindent=8mm \vspace{.1in}
Now we proceed to prove some useful properties of the operators $\mathcal{P}_{\alpha,j}$ and $\widetilde{\mathcal{P}}_{\alpha,j}$.

\parindent=0mm \vspace{.1in}
\textbf{Lemma 4.8.} The operators $\mathcal{P}_{\alpha,j}$ and $\widetilde{\mathcal{P}}_{\alpha,j}$ satisfy the following properties

\parindent=0mm \vspace{.1in}
(a) $\mathcal{P}_{\alpha,j} f = f $ if and only if $ f \in V_j^\alpha$ and $\widetilde{\mathcal{P}}_{\alpha,j} = f$ if only if $f \in \widetilde{V}_j^\alpha;$

\parindent=0mm \vspace{.1in}
(b) $\lim_{j\rightarrow \infty}\left\|\mathcal{P}_{\alpha,j}f -f\right\|_2 =0$ and $\lim_{j\rightarrow -\infty}\left\|\mathcal{P}_{\alpha,j}f\right\|_2 = 0$ for every $ f \in L^2(\mathbb R).$

\parindent=0mm \vspace{.1in}
\textbf{Proof.} (a) $\mathcal{P}_{\alpha,j} f = f$ if only if $f = \sum_{n\in \mathbb Z}\left\langle f ,\widetilde{\phi}_{\alpha, j,n}\right\rangle \phi_{\alpha,j,n}.$  Since $\left\{\phi_{\alpha,j,n} :n \in \mathbb Z\right\}$ is a Riesz basis for $V_j^\alpha$ and $\left\{\widetilde{\phi}_{\alpha,j,n} :n \in \mathbb Z\right\}$ is biorthogonal to $\{\phi_{\alpha,j,n} :n \in \mathbb Z\}$. By Lemma 3.8,  $f = \sum_{n\in \mathbb Z}\left\langle f ,\widetilde{\phi}_{\alpha,j,n}\right\rangle\phi_{\alpha,j,n}$ if and only if $f \in \overline{span}\{\phi_{\alpha,j,n} :n \in \mathbb Z\} = V_j^\alpha.$ Similar argument applies for $\widetilde{\mathcal{P}}_{\alpha,j}f$.

\parindent=0mm \vspace{.1in}
(b) It is straight forward.

\parindent=0mm \vspace{.1in}
\textbf{5. Biorthogonality of fractional wavelets}

\parindent=0mm \vspace{.1in}
\textbf{Definition 5.1.} Let $\phi$ and $\widetilde{\phi}$ be the scaling functions for dual fractional MRAs. Define the fractional wavelet $\psi$ and and $\widetilde{\psi}$ by 
 
\parindent=0mm \vspace{.0in}
$$ \psi(t) = \sum_{n \in \mathbb Z} g[n]\sqrt{2}\phi(2t-n)e^{\frac{-j}{2}[t^2 -(n/2)^2-(2t-n)^2]\cot\alpha}\eqno(5.1)$$
and
$$ \widetilde{\psi}(t) = \sum_{n \in \mathbb Z} \widetilde{g}[n]\sqrt{2}\widetilde{\phi}(2t-n)e^{\frac{-j}{2}[t^2 -(n/2)^2-(2t-n)^2]\cot\alpha}.\eqno(5.2)$$

\parindent=0mm \vspace{.0in}
The following lemma contains some basic properties of the wavelet and its dual.

\parindent=0mm \vspace{.1in}
\textbf{Lemma 5.2} Let $\psi$ and $\widetilde{\psi}$ be the wavelet and dual wavelet corresponding to the fractional MRA's $\{V_j^\alpha\}$ with the scaling function $\phi$ and $\{\widetilde{V}_j^\alpha\}$ with the scaling function $\widetilde{\phi}$. Then the following hold

\parindent=0mm \vspace{.1in}
(a) $\psi \in V_1^\alpha$ and $\widetilde{\phi} \in \widetilde{V}_1^\alpha$.

\parindent=0mm \vspace{.1in}
(b) $\left\{\widetilde{\psi}_{\alpha, 0,n} : n\in \mathbb Z\right\}$ is bi orthogonal to $\left\{\psi_{\alpha,0,n} : n\in \mathbb Z\right\}$.

\parindent=0mm \vspace{.1in}
(c) For all $m,n \in \mathbb Z$, we have 
 
\parindent=0mm \vspace{.1in}
$$\left \langle \psi_{\alpha,0,n} ,\widetilde{\phi}_{\alpha,0,m} \right\rangle = \left\langle \widetilde{\psi}_{\alpha,0,n}, \phi_{\alpha, 0,m}\right\rangle = 0.$$

\parindent=0mm \vspace{.1in}
\textbf{Proof.}(a) This clearly follows from the definition of $\psi$ and $\widetilde{\psi}$.

\parindent=0mm \vspace{.1in}
(b) Taking the FrFT on both sides of $(5.1)$ and $(5.2)$ gives

$$\mathcal{F}_\alpha\{\psi(t)\}(u)= \Gamma_\alpha\left(\frac{u}{2}\right)\Theta_\alpha\left(\frac{u}{2}\right);\eqno(5.3)$$

\parindent=0mm \vspace{.0in}
$$\mathcal{F}_\alpha\{\widetilde{\psi}(t)\}(u)=\widetilde{\Gamma}_\alpha\left(\frac{u}{2}\right)\widetilde{\Theta}_\alpha\left(\frac{u}{2}\right),\eqno(5.4)$$
where

$$\Gamma_\alpha(u) =e^{-2j\pi(u+2k\pi\sin\alpha)}\overline{\widetilde{\Lambda}_\alpha(u+2k\pi\sin\alpha)}; \eqno(5.5)$$

$$\widetilde{\Gamma}_\alpha(u) =e^{-2j\pi(u+2k\pi\sin\alpha)}\overline{\Lambda_\alpha(u+2k\pi\sin\alpha)}. \eqno(5.6)$$
Since $\{\phi_{\alpha,0,n} :n\in \mathbb Z\}$ is biorthogonal to $\{\widetilde{\phi}_{\alpha,0,n} :n\in \mathbb Z\}$,we have
\begin{align*}
\dfrac{1}{\sin\alpha} &= \sum_{n \in \mathbb Z}\Theta_\alpha(u+2n\pi\sin\alpha)\overline{\widetilde{\Theta}_\alpha(u+2n\pi\sin\alpha)}\\\
&=\sum_{n \in \mathbb Z}\Lambda_\alpha\left(\frac{u}{2}+2k\pi\sin\alpha\right)\Theta_\alpha\left(\frac{u}{2}+2n\pi\sin\alpha\right)\\\
&\qquad\qquad\qquad\qquad\qquad\qquad\qquad\times\overline{\widetilde{\Lambda}_\alpha\left(\frac{u}{2}+2k\pi\sin\alpha\right)}\overline{\widetilde{\Theta}_\alpha\left(\frac{u}{2}+2n\pi\sin\alpha\right)}\\\\
&=\Lambda_\alpha\left(\frac{u}{2}\right)\overline{\widetilde{\Lambda}_\alpha\left(\frac{u}{2}\right)}\sum_{k \in \mathbb Z}\Theta_\alpha\left(\frac{u}{2}+2k\pi\sin\alpha\right)\overline{\widetilde{\Theta}_\alpha\left(\frac{u}{2}+2k\pi\sin\alpha\right)}\\\
&~~~~~~~~ +\Lambda_\alpha\left(\frac{u}{2}+2k\pi\sin\alpha\right)\overline{\widetilde{\Lambda}_\alpha\left(\frac{u}{2} +2k\pi\sin\alpha\right)}\sum_{k \in \mathbb Z}\Theta_\alpha\left(\frac{u}{2}+(2k+1)\pi\sin\alpha\right)\\\
&\qquad\qquad\qquad\qquad\qquad\qquad\qquad\qquad\qquad\qquad\quad\times\overline{\widetilde{\Theta}_\alpha\left(\frac{u}{2}+(2k+1)\pi\sin\alpha\right)}\\\\
&=\Lambda_\alpha\left(\frac{u}{2}\right)\overline{\widetilde{\Lambda}_\alpha\left(\frac{u}{2}\right)}\dfrac{1}{\sin\alpha} +\Lambda_\alpha\left(\frac{u}{2}+2k\pi\sin\alpha\right)\overline{\widetilde{\Lambda}_\alpha\left(\frac{u}{2} +2k\pi\sin\alpha\right)}\dfrac{1}{\sin\alpha}.\tag{4.7}
\end{align*}

\parindent=0mm \vspace{.0in}
Combining (5.5), (5.6) and (5.7) gives
$$\Gamma_\alpha\left(\frac{u}{2}\right)\overline{\widetilde{\Gamma}_\alpha\left(\frac{u}{2}\right)}\dfrac{1}{\sin\alpha} +\Gamma_\alpha\left(\frac{u}{2}+2k\pi\sin\alpha\right)\overline{\widetilde{\Gamma}_\alpha\left(\frac{u}{2} +2k\pi\sin\alpha\right)}\dfrac{1}{\sin\alpha}=1.\eqno(5.8)$$
Mimicking the argument giving (4.7) gives

\begin{align*}
 &\sum_{n \in \mathbb Z}\mathcal{F}_\alpha\{\psi\}(u+2n\pi\sin\alpha)\overline{\mathcal{F}_\alpha\{\widetilde{\psi}\}(u+2n\pi\sin\alpha)} \qquad\qquad\qquad\qquad\qquad\qquad\\\
&\qquad =\Gamma_\alpha\left(\frac{u}{2}\right)\overline{\widetilde{\Gamma}_\alpha\left(\frac{u}{2}\right)}\,\dfrac{1}{\sin\alpha} +\Gamma_\alpha\left(\frac{u}{2} +2k\pi\sin\alpha\right)\overline{\widetilde{\Gamma}_\alpha\left(\frac{u}{2} +2k\pi\sin\alpha\right)}\dfrac{1}{\sin\alpha}\\\
 &\qquad=\dfrac{1}{\sin\alpha}.\tag{5.9}
\end{align*}
Therefore by Lemma 3.5.,$\left\{\psi_{\alpha,0,n} :n\in\mathbb Z\right\}$ is biorthogonal to $\left\{\widetilde{\psi}_{\alpha,0,n} :n\in\mathbb Z\right\}$.

\parindent=0mm \vspace{.1in}
(c) For fixed $m, n \in \mathbb Z$, we have by Plancherel's formula, we have
\begin{align*}
\left\langle \psi_{\alpha ,0,n} ,\widetilde{\phi}_{\alpha,0,m}\right\rangle &= \int_{-\infty}^{\infty} \mathcal{F}_\alpha\{\psi\}(u) \overline{\widetilde{\Theta}_\alpha(u)}e^{\frac{j}{2}(n^2-m^2)\cot\alpha-j(n-m)u\csc\alpha}du\\\\
&=\int_{-\infty}^{\infty} \Gamma_\alpha\left(\frac{u}{2}\right)\overline{\widetilde{\Lambda}_\alpha\left(\frac{u}{2}\right)}\Theta_\alpha\left(\frac{u}{2}\right)\overline{\widetilde{\Theta}_\alpha\left(\frac{u}{2}\right)}e^{\frac{j}{2}(n^2-m^2)\cot\alpha-j(n-m)u\csc\alpha}\,du\\\\
&=\dfrac{1}{\sin\alpha}\int_{0}^{2\pi\sin\alpha}e^{\frac{j}{2}(n^2 - m^2)\cot\alpha-j(n-m)u \csc\alpha}\sum_{k \in \mathbb Z}\Gamma_\alpha\left(\frac{u}{2}+ 2k\pi\sin\alpha\right)\\\
&\qquad\qquad\times\overline{\widetilde{\Lambda}_\alpha\left(\frac{u}{2} + 2k\pi\sin\alpha\right)}\Theta_\alpha\left(\frac{u}{2} + 2k\pi\sin\alpha\right)\overline{\widetilde{\Theta}_\alpha\left(\frac{u}{2} + 2k\pi\sin\alpha\right)}\,du \\\\
&=\dfrac{1}{\sin\alpha}\int_{0}^{2\pi\sin\alpha}e^{\frac{j}{2}(n^2 - m^2)\cot\alpha-j(n-m)u \csc\alpha}\Bigg\{ \Gamma_\alpha\left(\frac{u}{2}\right)\overline{\widetilde{\Lambda}_\alpha\left(\frac{u}{2}\right)} \\\
&~~~~~~~~~~~~~~~~~~~~~~~~~~~~~~~~~~~~~\times\sum_{\ell \in \mathbb Z}\Theta_\alpha\left(\frac{u}{2} + 2\ell\pi\sin\alpha\right)\overline{\widetilde{\Theta}_\alpha\left(\frac{u}{2} + 2\ell\pi\sin\alpha\right)}\\\
&~~~~~~+\Gamma_\alpha\left(\frac{u}{2}+2\ell\pi\sin\alpha\right)\overline{\widetilde{\Lambda}_\alpha\left(\frac{u}{2}+2\ell\pi\sin\alpha\right)}\\\
&~~~~~~~~~~~~~~~~~\times\sum_{\ell \in \mathbb Z}\Theta_\alpha\left(\frac{u}{2} + (2\ell+1)\pi\sin\alpha\right)\overline{\widetilde{\Theta}_\alpha\left(\frac{u}{2} + (2\ell+1)\pi\sin\alpha\right)}\Bigg\}du \\\\
&=\dfrac{1}{\sin^2\alpha}\int_{0}^{2\pi\sin\alpha}e^{\frac{j}{2}(n^2 - m^2)\cot\alpha-j(n-m)u \csc\alpha}\bigg\{\Gamma_\alpha\left(\frac{u}{2}\right)\overline{\widetilde{\Lambda}_\alpha\left(\frac{u}{2}\right)}\\\
&~~~~~~~~~~~~~~~~~~~~~~~~~~~~~~~~~~~~~~+\Gamma_\alpha\left(\frac{u}{2}+2\ell\pi\sin\alpha\right)\overline{\widetilde{\Lambda}_\alpha\left(\frac{u}{2}+2\ell\pi\sin\alpha\right)}\bigg\}\,du\\\
&=0.
\end{align*}
Similarly, we can show that $\left\langle \widetilde{\psi}_{\alpha,0,n},\phi_{\alpha,0,n}\right\rangle =0$ for all $n,m \in \mathbb Z.$

\parindent=8mm \vspace{.1in}
Now our main objective is to show that the wavelets associated with the dual fractional MRAs are biorthogonal and also they form Riesz basis for $L^2(\mathbb R).$ For that the following proposition is very useful.

\parindent=0mm \vspace{.1in}
\textbf{Proposition 5.3.} Let $\phi$ and $\widetilde{\phi}$ be the scaling functions for  dual fractional MRAs and $\psi,\widetilde{\psi}$ be the associated wavelets satisfying the matrix condition 
$$M_\alpha(u)\overline{\widetilde{M}_\alpha(u)} = I\eqno(5.10)$$
where
$$M_\alpha(u)=\begin{bmatrix} \Lambda_\alpha(u) &\Lambda_\alpha(u+2\pi\sin\alpha)\\ \Gamma_\alpha(u)&\Gamma_\alpha(u+2\pi\sin\alpha \end{bmatrix}\eqno(5.11)$$
 
Denote $\psi_0 = \phi$ and $\widetilde{\psi}_0 = \widetilde{\phi}$. Then for every $f \in L^2(\mathbb R)$, we have
$$\mathcal{P}_{\alpha,1}f = \mathcal{P}_{\alpha,0}f + \sum_{k \in \mathbb Z} \left\langle f, \widetilde{\psi}_{\alpha,0,k}\right\rangle \psi_{\alpha, 0,k}\eqno(5.12)$$
and
$$\widetilde{\mathcal{P}}_{\alpha,1}f = \widetilde{\mathcal{P}}_{\alpha,0}f + \sum_{k \in \mathbb Z} \left\langle f, \psi_{\alpha,0,k}\right\rangle \widetilde{\psi}_{\alpha, 0,k}\eqno(5.13)$$
where the series converges in $L^2(\mathbb R)$

\parindent=0mm \vspace{.1in}
\textbf{Proof.} We will prove only (5.12) as the proof of (5.13) follows in the similar manner. Further it suffices to prove (5.12) in the weak sense, that is, for all $f,g \in L^2(\mathbb R)$,
\begin{align*}
\left\langle \mathcal{P}_{\alpha,1}f,g \right\rangle & = \left\langle \mathcal{P}_{\alpha,0}f,g \right\rangle +\sum_{k \in \mathbb Z} \left\langle f, \widetilde{\psi}_{\alpha,0,k}\right\rangle \overline{\left\langle g, \psi_{\alpha,0,k}\right\rangle} \\\
&= \sum_{k \in \mathbb Z} \left\langle f, \widetilde{\psi}_{\alpha,0,k}\right\rangle \overline{\left\langle g, \psi_{\alpha,0,k}\right\rangle}. 
\end{align*}
We have
\begin{align*}
&\sum_{k \in \mathbb Z} \left\langle f, \widetilde{\psi}_{\alpha,0,k}\right\rangle \overline{\left\langle g, \psi_{\alpha,0,k}\right\rangle}\\\
 &\qquad\qquad=\sum_{k\in \mathbb Z}\left\{ \int_{-\infty}^{\infty}\mathcal{F}_\alpha\{f(t)\}(u) \overline{\mathcal{F}_{\alpha}\{\widetilde{\psi}\}(u)}e^{\frac{jn^2}{2}\cot\alpha-jnu\csc\alpha}\,du\right\} \\\
&\qquad\qquad\qquad\qquad\qquad\qquad\qquad\times\left\{ \int_{-\infty}^{\infty}\overline{\mathcal{F}_\alpha\{g(t)\}(u)} \overline{\mathcal{F}_{\alpha}\{\psi\}(u)}e^{\frac{-jn^2}{2}\cot\alpha+jnu\csc\alpha}\,du\right\}\\\\
&\qquad\qquad=\sum_{k\in \mathbb Z}\left\{ \int_{0}^{2\pi\sin\alpha}\sum_{\ell \in \mathbb Z}\mathcal{F}_\alpha\{f(t)\}(u+2\ell\pi\sin\alpha)\right. \\\
&\qquad\qquad\qquad\qquad\qquad\qquad\qquad\times\left.\overline{\mathcal{F}_{\alpha}\{\widetilde{\psi}\}(u+2\ell\pi\sin\alpha)}e^{\frac{jn^2}{2}\cot\alpha-jn(u+2\ell\pi\sin\alpha)\csc\alpha}\,du\right\} \\\
&\qquad\qquad\qquad\qquad\times\left\{ \int_{0}^{2\pi\sin\alpha}\sum_{\ell'\in\mathbb Z}\overline{\mathcal{F}_\alpha\{g(t)\}(u+2\ell'\pi\sin\alpha)}\right. \\\
&\qquad\qquad\qquad\qquad\qquad\qquad\quad\times\left.\mathcal{F}_{\alpha}\{\psi\}(u+2\ell'\pi\sin\alpha)e^{\frac{-jn^2}{2}\cot\alpha+jn(u+2\ell'\psi\sin\alpha)\csc\alpha}\,du\right\}\\\\
&\qquad\qquad=\sum_{k\in \mathbb Z}\left\{ \int_{0}^{2\pi\sin\alpha}\sum_{\ell \in \mathbb Z}\mathcal{F}_\alpha\{f(t)\}(u+2\ell\pi\sin\alpha) \overline{\mathcal{F}_{\alpha}\{\widetilde{\psi}\}(u+2\ell\pi\sin\alpha)}\,du\right\} \\\
&\qquad\qquad\qquad\qquad\times\left\{ \int_{0}^{2\pi\sin\alpha}\sum_{\ell'\in\mathbb Z}\overline{\mathcal{F}_\alpha\{g(t)\}(u+2\ell'\pi\sin\alpha)} \mathcal{F}_{\alpha}\{\psi\}(u+2\ell'\pi\sin\alpha)\,du\right\}\\\\
&\qquad\qquad=\int_{0}^{2\pi\sin\alpha}\sum_{\ell \in \mathbb Z}\mathcal{F}_\alpha\{f(t)\}(u+2\ell\pi\sin\alpha)\overline{\widetilde{\Gamma}_\alpha\left(\frac{u}{2}+2\ell\pi\sin\alpha\right)}\overline{\widetilde{\Theta}_\alpha\left(\frac{u}{2} +2 \ell\pi\sin\alpha\right)}\\\
&\qquad\qquad\qquad\quad\times\sum_{\ell' \in \mathbb Z}\overline{\mathcal{F}_\alpha\{g(t)\}(u+2\ell'\pi\sin\alpha)}\widetilde{\Gamma}_\alpha\left(\frac{u}{2}+2\ell\pi\sin\alpha\right)\Theta_\alpha\left(\frac{u}{2} +2 \ell\pi\sin\alpha\right)du\\\\
&\qquad\qquad=\int_{0}^{2\pi\sin\alpha}\sum_{\ell \in \mathbb Z}\sum_{\ell' \in \mathbb Z}\mathcal{F}_\alpha\{f\}(u+2\ell\pi\sin\alpha)\overline{\widetilde{\Theta}_\alpha\left(\frac{u}{2}+2\ell\sin\alpha\right)}\\\
&\qquad\qquad\qquad\qquad\qquad\qquad\quad\times\overline{\mathcal{F}_\alpha\{g\}(u+2\ell'\pi\sin\alpha)}\Theta_\alpha\left(\frac{u}{2}+2\ell'\pi\sin\alpha\right)\,du.\tag{5.14}
\end{align*}

In the similar lines, we can obtain
\begin{align*}
\sum_{k\in \mathbb Z}\left\langle f, \widetilde{\phi}_{\alpha, 1,k}\right\rangle\overline{\left\langle g,\phi_{\alpha,1,k}\right\rangle}& = \int_{0}^{2\pi\sin\alpha}\sum_{s\in \mathbb Z}\sum_{s'\in \mathbb Z} \mathcal{F}_\alpha\{f\}(u+2s\pi\sin\alpha),\overline{\widetilde{\Theta}_\alpha\left(\frac{u}{2}+ 2s\pi\sin\alpha\right)}\\\
&\qquad\quad\times\overline{\mathcal{F}_\alpha\{g\}(u+2s'\pi\sin\alpha)}\Theta_\alpha\left(\frac{u}{2}+2s'\pi\sin\alpha\right)\,du.\tag{5.15}
\end{align*}
Since the right hand sides of $(5.14)$ and $(5.15)$ are same. This completes the proof.$\square$

\parindent=8mm \vspace{.1in}
Combining the Proposition 5.3. and Lemma 4.8., we have the following proposition.

\parindent=0mm \vspace{.1in}
\textbf{Proposition 5.4.} Let $\phi,\widetilde{\phi}$ and $\psi,\widetilde{\psi}$ be defined as above. Then for every $f \in L^2(\mathbb R),$ we have
 
$$ f = \sum_{j \in \mathbb Z}\sum_{j \in \mathbb Z}\left\langle f, \widetilde{\psi}_{\alpha,j,k}\right\rangle \psi_{\alpha,j,k} =\sum_{j \in \mathbb Z}\sum_{j \in \mathbb Z}\left\langle f, \psi_{\alpha,j,k}\right\rangle \widetilde{\psi}_{\alpha,j,k}\eqno(5.16)$$

\textbf{Theorem 5.5.} Let $\phi$ and $\widetilde{\phi}$ be the scaling functions for Dual fractional MRAs and $\psi,\widetilde{\psi}$ be the associated wavelets as in the Proposition 4.3. Then the collections $\{\psi_{\alpha,j,k} :j,k \in \mathbb Z\}$ and $\{\widetilde{\psi}_{\alpha,j,k} : j,k \in \mathbb Z\}$ are biorthogonal. Further, if 
 
$$ |\Theta_\alpha(u)| \le C\left(1+|u|\right)^{-1/2-\epsilon},~~~~~~~|\widetilde{\Theta}_\alpha(u)| \le C \left(1+|u|\right)^{-1/2-\epsilon},\eqno(5.17)$$
  
$$\left|\mathcal{F}_\alpha\{\psi\}(u)\right| \le C|u|,~~\textit{and}~~\left|\mathcal{F}_\alpha\{\widetilde{\psi}\}(u)\right| \le C|u|\eqno(5.18)$$
for some constant $C >0, \epsilon >0$ and for a.e., $u \in \mathbb R,$ then the collections $\{\psi_{\alpha,j,k} :j,k \in \mathbb Z\}$ and $\{\widetilde{\psi}_{\alpha,j,k} : j,k \in \mathbb Z\}$ form Riesz basis for $L^2(\mathbb R)$.

\parindent=0mm \vspace{.1in}
\textbf{Proof.} We begin the proof by proving that the collections $\left\{\psi_{\alpha,j,k} :j,k \in \mathbb Z\right\}$ and $\left\{\widetilde{\psi}_{\alpha,j,k} : j,k \in \mathbb Z\right\}$ are biorthogonal to each other. First we will show that, for $j \in \mathbb Z$
$$\left\langle \psi_{\alpha, j,k}, \widetilde{\psi}_{\alpha,j,k'}\right\rangle = \delta_{k,k'}.$$
We have already proved it for $j=0$, by Lemma 5.2. (b). For $j\ne 0$, we have
\begin{align*}
\left\langle \psi_{\alpha, j,k}, \widetilde{\psi}_{\alpha,j,k'}\right\rangle &=\left\langle \delta_{-j}\psi_{\alpha, 0,k}, \delta_{-j}\widetilde{\psi}_{\alpha,0,k'}\right\rangle \\\
& =\left\langle \psi_{\alpha, 0,k}, \widetilde{\psi}_{\alpha,0,k'}\right\rangle \\\
& = \delta_{k,k'}.
\end{align*}
Let $k,k' \in \mathbb Z$ be fixed and Let $j,j' \in \mathbb Z$. Assume that $j < j'$, we will show that

$$\left\langle \psi_{\alpha, j,k}, \widetilde{\psi}_{\alpha,j',k'}\right\rangle=0.$$

It can be shown that $\psi_{\alpha,0,k} \in V_1^\alpha$. Hence,$\psi_{\alpha,j,k} = \delta_{-j}\psi_{\alpha,0,k} \in V_{j+1}^\alpha \subseteq V_{j'}^\alpha$. Therefore, it will be enough to show that $\widetilde{\psi}_{\alpha,j',k'}$ is orthogonal to every element of $V_{j'}^\alpha$. Let $f \in V_{j'}^\alpha$. By Lemma 4.3, $\{\phi_{\alpha,j',k'} : k\in \mathbb Z\}$ is a Riesz basis for $V_{j'}^\alpha$. Hence, there exists a sequence $c[k] \in \ell^2(\mathbb Z)$ such that
$$ f = \sum_{k \in \mathbb Z} c[k] \phi_{\alpha, j',k'}~~~\textit{in}~~L^2(\mathbb R).$$

By Lemma 5.2 (c), we have 
\begin{align*}
\left\langle \widetilde{\psi}_{\alpha,j',k'}, \phi_{\alpha,j',k}\right\rangle &= \left\langle \delta_{-j'}\widetilde{\psi}_{\alpha,0,k'}, \delta_{-j'}\phi_{\alpha,0,k}\right\rangle\\\
& =\left\langle \widetilde{\psi}_{\alpha,0,k'}, \phi_{\alpha,0,k}\right\rangle\\\
&=0.
\end{align*}
Hence,
\begin{align*}
\left\langle \widetilde{\psi}_{\alpha,j',k'}, f\right\rangle &=\left\langle \widetilde{\psi}_{\alpha,j',k'},\sum_{k \in \mathbb Z} c[k] \phi_{\alpha, j',k'} \right\rangle\\\
&=\sum_{k \in \mathbb Z}\overline{c[k]}\left\langle \widetilde{\psi}_{\alpha,j',k'}, \phi_{\alpha,j',k}\right\rangle\\
&=0.
\end{align*}
In order to show that these two collections form Riesz bases for $L^2(\mathbb R)$, we must verify that they are linearly independent and satisfy the frame condition. Since they are biorthogonal to each other,therefore by Lemma 3.3, both the collections are linearly independent.

\parindent=8mm \vspace{.1in}
In order to show the frame condition, we must show that there exist constants $A, B, \widetilde{A},$ and $\widetilde{B} >0$ such that for every $f \in L^2(\mathbb R),$ we have

$$A\|f\|_2^2 \le \sum_{j\in \mathbb Z} \sum_{k \in \mathbb Z} \left|\left\langle f, \psi_{\alpha,j,k}\right\rangle \right|^2 \le B\|f\|_2^2,\eqno(5.18),$$
and
$$\widetilde{A}\|f\|_2^2 \le \sum_{j\in \mathbb Z} \sum_{k \in \mathbb Z} \left|\left\langle f, \widetilde{\psi}_{\alpha,j,k}\right\rangle \right|^2 \le \widetilde{B}\|f\|_2^2,\eqno(5.19).$$
We first establish the existence of upper bounds in (5.18) and (5.19). we have
\begin{align*}
\sum_{k \in \mathbb Z} \left|\left\langle f, \psi_{\alpha,j,k}\right\rangle\right|^2& =\sum_{k \in \mathbb Z} \left|\int_{-\infty}^{\infty} \mathcal{F}_\alpha\{f\}(u)\overline{\mathcal{F}_\alpha\{\psi\}(2^{j}u)} \,du\right|^2 \\\
&=\sum_{k \in \mathbb Z} \left|\int_{0}^{2\pi\sin\alpha}\sum_{m\in\mathbb Z} \mathcal{F}_\alpha\{f\}(u+2^j 2\pi m\sin\alpha)\overline{\mathcal{F}_\alpha\{\psi\}(2^{j}u+2m\pi\sin\alpha)} \,du\right|^2 \\\
&\le \int_{-\infty}^{\infty}\left|\mathcal{F}_\alpha\{f\}(u)\right|^2\,\left|\mathcal{F}_\alpha\{\psi\}(2^{j}u)\right|^{2\delta}\sum_{n\in \mathbb Z}\left|\mathcal{F}_\alpha\{\psi\}(2^{j}u+2n\pi\sin\alpha)\right|^{2(1-\delta)}\,du.
\end{align*}
We have assumed that $\left|\Theta_\alpha(u)\right| \le C\left(1+|u|\right)^{\frac{-1}{2}-\epsilon}$, hence we have $\left|\mathcal{F}_\alpha\{\psi\}(u)\right| \le C \left(1+|\frac{u}{2}|\right)^{\frac{-1}{2}-\epsilon}$. Therefore, $\sum_{n \in \mathbb Z}\left|\mathcal{F}_\alpha\{\psi\}(2^{j}u+2n\pi\sin\alpha)\right|^{2(1-\delta)}$ is uniformly bounded if $\delta < 2\epsilon(1+2\epsilon)^{-1}$. Hence, there exists $C>0$ such that
\begin{align*}
\sum_{j\in \mathbb Z}\sum_{k \in \mathbb Z} \left|\left\langle f, \psi_{\alpha,j,k}\right\rangle\right|^2& \le C\int_{-\infty}^{\infty}\left|\mathcal{F}_\alpha\{f\}(u)\right|^2\,\left|\mathcal{F}_\alpha\{\psi\}(2^{j}u)\right|^{2\delta}\,du\\\
& \le C\sup\left\{\sum_{j\in\mathbb Z}\left|\mathcal{F}_\alpha\{\psi\}(2^{j}u)\right|^{2\delta}: u \in [0, 2\pi\sin\alpha] \right\}\|f\|_2^2.
\end{align*}
Further for $1 < |u| \le 2\sin\alpha$, we have
\begin{align*}
\sum_{j=-\infty}^{0}\left|\mathcal{F}_\alpha\{\psi\}(2^{j}u)\right|^{2\delta}&\le \sum_{j=0}^{\infty}\frac{C^{2\delta}}{(1+|2^{j-1}u|)^{\delta(1+2\epsilon)}}\\\
&\le\sum_{j=0}^{\infty}\frac{C^{2\delta}}{2^{(j-1)\delta(1+2\epsilon)}}\\\
& =C^{2\delta}\frac{2^{\delta(1+2\epsilon)}}{1-2^{-\delta(1+2\epsilon)}}.
\end{align*}
 Also
\begin{align*}
\sum_{j=1}^{\infty}\left|\mathcal{F}_\alpha\{\psi\}(2^{j}u)\right|^{2\delta} &\le \sum_{j=1}^{\infty}\left(C2^{-j}|u|\right)^{2\delta}\\\
&\le C^{2\delta}\sum_{j=1}^{\infty}2^{(-j+1)2\delta}\\\
& =C^{2\delta}\frac{1}{1-2^{-2\delta}} 
\end{align*}
These two estimates show that $\sup\left\{\sum_{j\in\mathbb Z}\left|\mathcal{F}_\alpha\{\psi\}(2^{j}u)\right|^{2\delta}: u \in [0, 2\pi\sin\alpha] \right\}$ is finite. Hence, there exists $B>0$ such that the second inequality in (5.18) holds. In the similar manner we can obtain the upper bound in (5.19).

\parindent=8mm \vspace{.1in}
Now we proceed to show the existence of the lower bounds in (5.18) and (5.19), by virtue of the existence of upper bounds. From the proposition 5.4, if $f \in L^2(\mathbb R)$, then we have
$$ f = \sum_{j \in \mathbb Z}\sum_{k \in \mathbb Z}\left\langle f, \widetilde{\psi}_{\alpha,j,k}\right\rangle \psi_{\alpha,j,k} =\sum_{j \in \mathbb Z}\sum_{k \in \mathbb Z}\left\langle f, \psi_{\alpha,j,k}\right\rangle \widetilde{\psi}_{\alpha,j,k}.$$
Therefore, we have
\begin{align*}
\left\|f\right\|^2 & = \left\langle f, f\right\rangle \\\
&=\left\langle\sum_{j \in \mathbb Z}\sum_{k \in \mathbb Z}\left\langle f, \widetilde{\psi}_{\alpha,j,k}\right\rangle \psi_{\alpha,j,k}, f\right\rangle \\\
&=\sum_{j \in \mathbb Z}\sum_{k \in \mathbb Z}\left\langle f, \widetilde{\psi}_{\alpha,j,k}\right\rangle \left\langle\psi_{\alpha,j,k}, f\right\rangle \\\
& \le \left\{\sum_{j \in \mathbb Z}\sum_{k \in \mathbb Z}\left|\left\langle f, \widetilde{\psi}_{\alpha,j,k}\right\rangle\right|^2\right\}^{1/2}\left\{\sum_{j \in \mathbb Z}\sum_{k \in \mathbb Z}\left|\left\langle f, \psi_{\alpha,j,k}\right\rangle\right|^2\right\}^{1/2}\\\
&\le\sqrt{\widetilde{B}}\|f\|_2\left\{\sum_{j \in \mathbb Z}\sum_{k \in \mathbb Z}\left|\left\langle f, \psi_{\alpha,j,k}\right\rangle\right|^2\right\}^{1/2}.
\end{align*}
Hence,

$$\dfrac{1}{\widetilde{B}}\|f\|_2^2 \le \sum_{j \in \mathbb Z}\sum_{k \in \mathbb Z}\left|\left\langle f, \psi_{\alpha,j,k}\right\rangle\right|^2.$$
Similarly we can show that
$$\dfrac{1}{B}\|f\|_2^2 \le \sum_{j \in \mathbb Z}\sum_{k \in \mathbb Z}\left|\left\langle f, \widetilde{\psi}_{\alpha,j,k}\right\rangle\right|^2.$$
This completes the proof of the theorem.

\parindent=0mm \vspace{.1in}

{\bf{References}}

\end{document}